\documentclass[11pt]{amsart}
\usepackage[a4paper,left=2cm,right=2cm,bottom=3cm]{geometry}

\usepackage[latin1]{inputenc}
\usepackage[english]{babel}

\usepackage{amsmath,amssymb}
\usepackage{amsthm}
\usepackage{tikz}

\newtheorem{theorem}{Theorem}[]

\begin{document}

\title{The reverse isoperimetric inequality for convex shapes included in a triangle}
\author[]{Rida Ait El Manssour}
\author[]{Zakaria Fattah}
\author[]{Ilias Ftouhi}

\address[Rida Ait El Manssour]{}
\email{rida.aitelmanssour@cs.ox.ac.uk}

\address[Zakaria Fattah]{Laboratoire MIPA, N\^imes Universit\'e, N\^imes, France.}
\email{zakaria.fattah@unimes.fr}

\address[Ilias Ftouhi]{Laboratoire MIPA, N\^imes Universit\'e, N\^imes, France.}
\email{ilias.ftouhi@unimes.fr}

\date{\today}

\begin{abstract}
We establish a sharp reverse isoperimetric inequality for convex
subsets of a triangle under a prescribed area constraint.
Given a triangle $T$ and an area fraction
$\eta\in(0,1)$, we determine the set with maximal perimeter among convex
bodies $K\subset T$ satisfying $|K|=\eta|T|$ and characterize
the equality cases. Every maximizer is a triangle, and the optimal
configuration is obtained by retaining a longest side of $T$
and placing the third vertex on a shortest side, at the position
determined by the area constraint. This characterization yields
an explicit expression for the maximal perimeter in terms of
$\eta$ and the side lengths of $T$. 
\end{abstract}



\maketitle


\section{Introduction}

The classical isoperimetric inequality in the plane asserts that
\[
P(K)^2\geq 4\pi |K|
\]
for every measurable planar set $K\subset \mathbb{R}^2$, with equality if and only if $K$ is
a disk. Here and throughout the paper, $|K|$ and $P(K)$ denote the
area and the perimeter of $K$, respectively. Thus, among sets of
prescribed area, the disk minimizes the perimeter; see
\cite{Blasjo2005} for a historical account of the isoperimetric
problem. The corresponding maximization problem is fundamentally
different as even within the class of convex sets, the perimeter is
unbounded under an area constraint alone.
Therefore, a potential reverse isoperimetric inequality cannot hold unless additional constraints are considered. 

The literature on reverse isoperimetric inequalities is quite rich. For example, in the affine setting, Ball \cite{Ball1991} proved that every convex body in $\mathbb{R}^n$ admits an affine image whose isoperimetric ratio is not larger than that of a regular simplex. For centrally symmetric convex bodies, the corresponding extremal body is the cube. These results concern a suitable representative of each affine equivalence class, rather than an upper bound valid for every convex body in its original position.

Another approach is to impose restrictions on curvature. Howard and Treibergs \cite{HowardTreibergs1995} established reverse
isoperimetric inequalities for plane curves with bounded curvature, together with stability results. More recently, Croce, Fattah and Pisante \cite{CroceFattahPisante2022} considered planar $(\alpha,\beta)$-convex bodies, whose radius of curvature is bounded between two prescribed values $\alpha$ and $\beta$. They characterized the extremal shapes called $(\alpha,\beta)$-eggs that are sets whose boundary consists of two circular arcs of radius $\alpha$ joined by two circular arcs
of radius $\beta$. Their work generalizes (at least in the planar case) a result by Borisenko and Drach \cite{zbMATH06434370} that concerns only an upper bound on the radius of curvature but in any dimensions. 

In the present paper, in addition to convexity, we assume an inclusion constraint. More precisely, given a planar and convex body $\Omega$ and a prescribed area $m\in(0,|\Omega|)$, we will consider problems of the type
\[
\max\bigl\{P(K): K\subset\Omega\text{ is a convex body},
\ |K|=m\bigr\}.
\]
The inclusion constraint prevents the elongation responsible for the unboundedness of the unconstrained problem. Moreover, Blaschke's selection theorem (see for example \cite[Thm.
1.8.7]{Schneider2014}) and the continuity of the area and perimeter functionals under the Hausdorff convergence of convex bodies imply the existence of an extremal shape. 

This problem has classical antecedents in Favard's work on extremal convex curves \cite{Favard1929}. For a disk as a container, Bogosel \cite{Bogosel2024} recently gave a complete characterization of the
optimal shapes: A maximizer is a polygon inscribed in the disk, with all but one of its sides equal and the remaining side no longer than the others. Its number of sides and its side lengths are determined by the prescribed area. In particular, the number of sides changes as the area crosses the areas of regular inscribed polygons. For a convex polygonal container with $N$ sides, Bogosel also remarks the existence of a polygonal maximizer with at most $2N$ vertices, all lying on the boundary of the container. Related inclusion-constrained problems were studied by Bianchini and Henrot \cite{BianchiniHenrot2012}, who were interested in the minimization of the functional  $\lambda|K|-P(K)$ over convex sets constrained between two concentric disks. We finally, note that a further connection is provided in the recent work by Fattah, Ftouhi and Zuazua \cite{FattahFtouhiZuazua2025}, where the authors study the problem of minimizing a so called $L^p$-distance between a convex container and a convex subset of a given fraction of area.  

Inspired by \cite{Bogosel2024}, we aim in the present paper to resolve the reverse isoperimetric problem in the case of a triangular container. More precisely, let
$T$ be a nondegenerate triangle. We are interested in the following problem 
\begin{equation}\label{problem}
  \max\{P(K) |\ K\ \text{convex}\subset T,\ \text{and}\ |K|=\eta |T|\}.  
\end{equation}

Our main result is the following:


\begin{theorem}
\label{thm:reverse-isoperimetric-triangle}
Let $T$ be a nondegenerate triangle of vertices $A$, $B$ and $C$ with side lengths
\[
a:=|BC|\geq b:=|CA|\geq c:=|AB|>0,
\]
and let $\eta\in(0,1)$. For every convex body $K\subset T$ satisfying $|K|=\eta|T|$, we have
\begin{equation*}
\label{eq:reverse-isoperimetric-triangle}
P(K) \leq a+\eta c+ \sqrt{(1-\eta)a^2+\eta b^2-\eta(1-\eta)c^2}.
\end{equation*}

The inequality is sharp and is attained only for the triangle $K_\eta$ of vertices $B$, $C$ and $D_\eta$, where $D_\eta:=(1-\eta)B+\eta A$, or any other triangle congruent to $K_\eta$ which is also included in $T$. 
\end{theorem}

\begin{figure}[h!]
    \centering
    \includegraphics[width=0.5\linewidth]{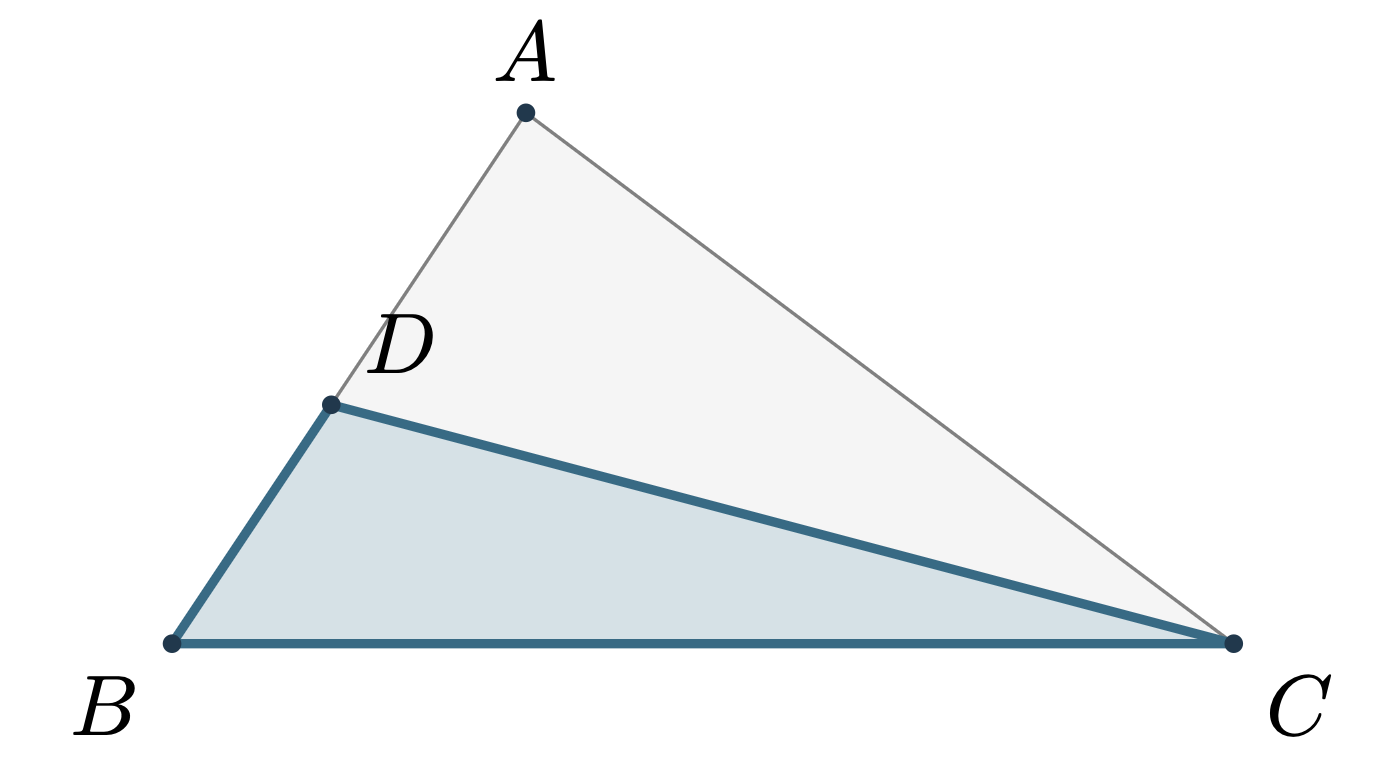}
    \caption{The set maximizing the perimeter among convex subsets of $T$ of a given area.}
    \label{fig:placeholder}
\end{figure}

The proof is presented in the next section and is based on the polygonal nature of solutions of \eqref{problem} proved in \cite{LamboleyNovruzi2010} (see also \cite{FattahFtouhiZuazua2025}) combined with relevant geometric transformations that preserve the area while increasing the perimeter such as \text{parallel chord movements}.

\section{Proof of the main theorem}
The existence of an optimal shape  is a direct consequence of Blaschke's selection theorem \cite[Thm. 1.8.7]{Schneider2014} and the continuity of the area and the perimeter for convex bodies under the Hausdorff convergence. We will present the proof in 6 steps. \newpage

\textbf{Step 1: The optimal shape is a polygon of at most 6 vertices, all located on $\partial T$} \vspace{2mm}

In \cite[Proposition 11]{FattahFtouhiZuazua2025} it has been shown that for any solution $K$ of problem \eqref{problem}, the set $\partial K\cap \mathring{T}$ is given by the union of disjoint segments. Since $T$ is a triangle and $K$ is convex, the set $\partial K\cap \mathring{T}$ is given by the union of at most $3$ segments. We finally deduce that an optimal shape $K$ should be a polygon with at most $6$ vertices, all located on the boundary $\partial T$.  \vspace{2mm}

\textbf{Step 2: Excluding hexagons and pentagons}\vspace{2mm}

\begin{figure}[h!]
    \centering
    \includegraphics[width=0.6\linewidth]{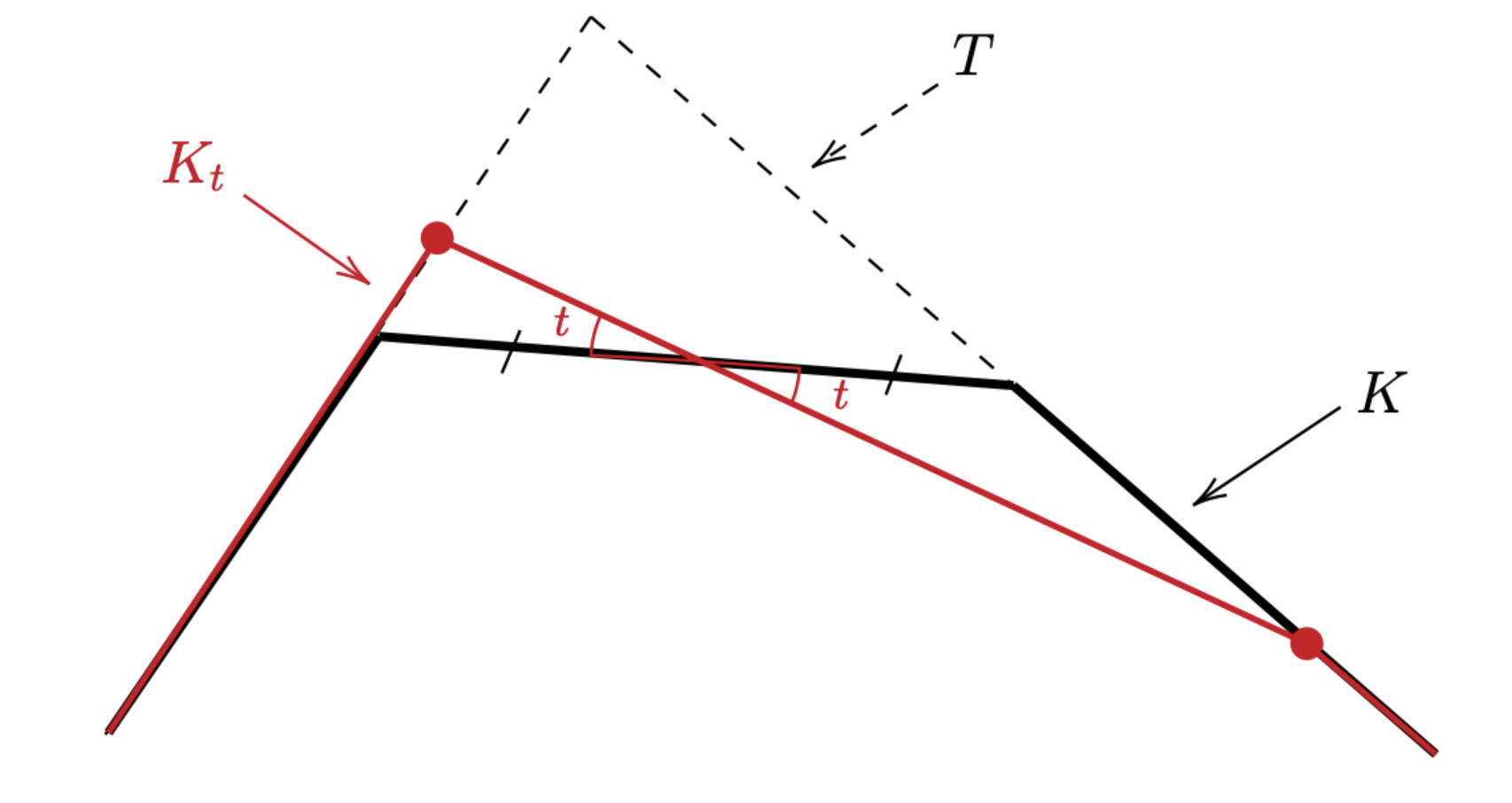}
    \caption{The polygon $K_t$ obtained by rotating the middle side of $K$ along its midpoint.}
    \label{fig:bucur}
\end{figure}

In this step, we prove that a domain $K$ which is a solution of Problem \eqref{problem} cannot have more than four sides. Indeed, if $K$ has $5$ or $6$ sides, then, there are three consecutive sides such as the first and the third are located on two consecutive sides of the triangle $T$ and the second contained in $\mathring{T}$. We consider a perturbation introduced in \cite{BF} which transforms $K$ into a new polygon $K_t$ obtained by rotating the middle side along its midpoint as shown in Figure \ref{fig:bucur}. For sufficiently small $t$, it has been shown in the second step of the proof of  \cite[Lemma 18]{Ftouhi2023} that for sufficiently small $t>0$, we have 
$$|K_t|\leq |K|\ \ \ \text{and}\ \ \ P(K_t)>P(K).$$

We can then construct a convex set $\Tilde{K_t}$ such that $|\Tilde{K_t}|= |K|= \eta |T|$ and $K_t\subset \Tilde{K_t}\subset T$. Thus, since the perimeter of convex sets is increasing with respect to inclusions, we deduce that 
$$|\Tilde{K_t}|= |K|\ \ \ \text{and}\ \ \ P(\Tilde{K_t})>P(K).$$

Therefore, a polygon $K$ with more than $4$ sides cannot be a maximizer in \eqref{problem} as we can always construct a better set. 

\textbf{Step 3: Excluding quadrilaterals}\vspace{2mm}

Let us now show that a quadrilateral $K$ cannot be a solution of Problem \eqref{problem}. If there are two sides of $T$ each containing two vertices of $K$, then, three consecutive sides of $K$ are in the same position presented in Figure \ref{fig:bucur}. We can therefore use the same argument above to show that $K$ cannot be a maximizer in \eqref{problem}. We then assume that each side of $T$ contains at least one vertex of $K$, which means one of the sides of $T$ contains two vertices of $K$, meanwhile, the interior of the two remaining sides contain exactly one vertex each. 

\begin{figure}[h!]
    \centering
    \includegraphics[width=0.5\linewidth]{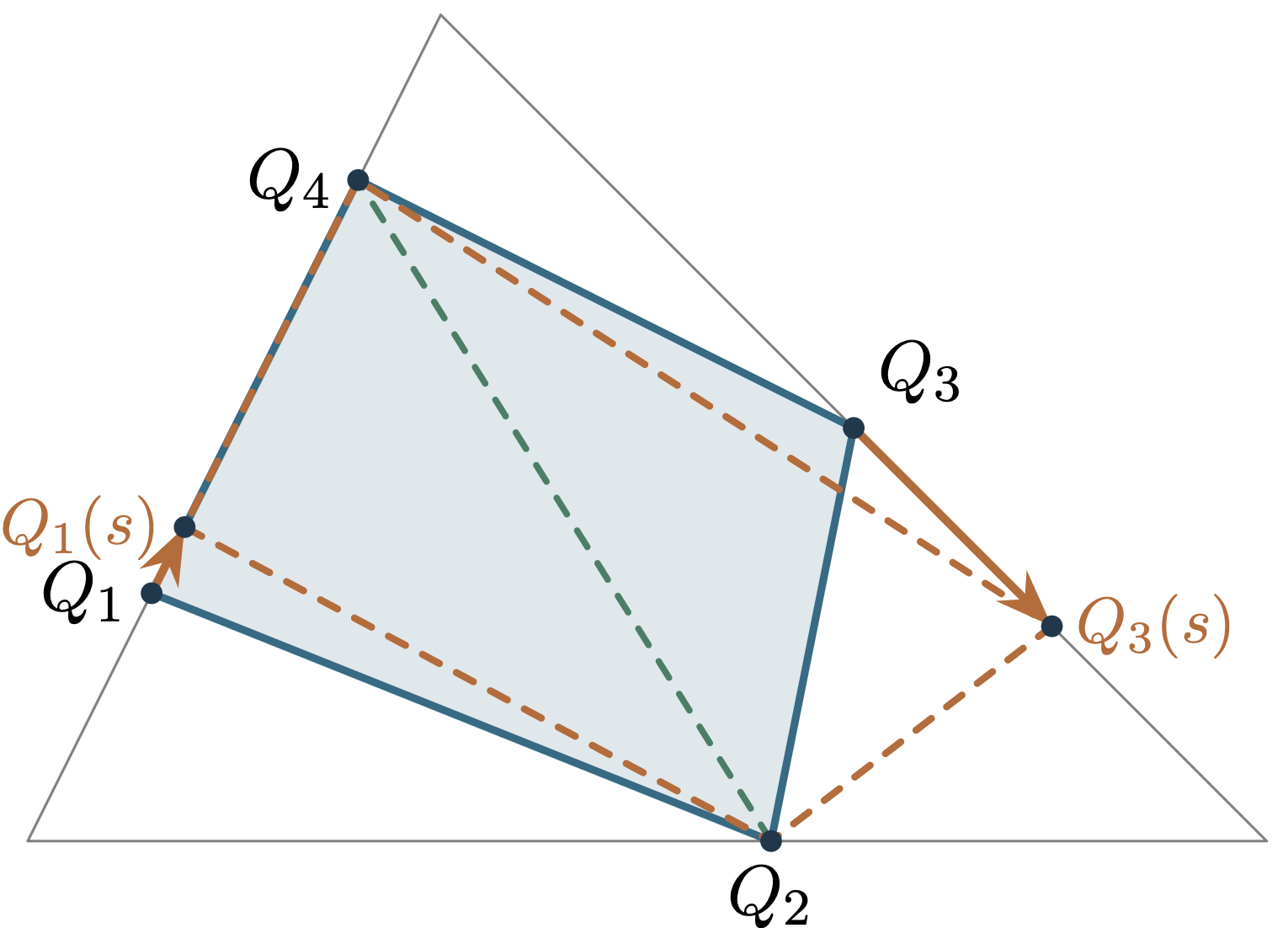}
    \caption{A quadrilateral and its perturbation where the area is preserved.}
    \label{fig:quadri}
\end{figure}

It remains then to exclude the case where the quadrilateral $K$ has two non adjacent vertices located in the interior of two different sides of $T$. We then denote by $Q_1$, $Q_2$, $Q_3$ and $Q_4$ the successive vertices of the convex $K$ and assume that the (nonadjacent) vertices $Q_2$ and $Q_4$ are in the interior of two sides of $T$, see Figure \ref{fig:quadri}. Let us consider two vectors $\Vec{u}$ and $\Vec{v}$ that are respectively collinear to $(Q_1Q_4)$ and the the side of $T$ containing $Q_3$ and define $\mu\in \mathbb{R}$. For sufficiently small $s$, we define $Q_1(s):=Q_1+s\Vec{u}$ and  $Q_3(s):=Q_3+s\mu \Vec{v}$, where $\mu\in \mathbb{R}$ is chosen such that $\Vec{u}-\mu \Vec{v}$ is collinear to the vector $\overrightarrow{Q_4Q_2}$. We then obtain a perturbed quadrilateral $K_s$ of vertices $Q_1(s)$, $Q_2$, $Q_3(s)$ and $Q_4$, such that for every $s$ sufficiently small,
\begin{align*}
|K_s|-|K| &= \frac{1}{2}\det(\overrightarrow{Q_3Q_1},\overrightarrow{Q_4Q_2}) - \frac{1}{2}\det(\overrightarrow{Q_3(s)Q_1(s)},\overrightarrow{Q_4Q_2}) \\
&= \frac{s}{2} \det(\Vec{u}-\mu \Vec{v},\overrightarrow{Q_4Q_2}) = 0. 
\end{align*}

Let us now see what happens to the perimeters. We write 
$$P(K_s)-P(K) = (Q_1(s)Q_2 + Q_2Q_3(s) + Q_3(s)Q_4 + Q_4Q_1(s)) - (Q_1Q_2 + Q_2Q_3 + Q_3Q_4 + Q_4Q_1).$$

The function $s\longmapsto P(K_s)-P(K)$ is convex as it is given by the sum of a constant and functions of the type $\Psi:s\longmapsto |x+sw|$, where $x$ and $w$ are given vectors of $\mathbb{R}^2$. For all sufficiently small $s$, we have 
$$\Psi''(s) = \frac{\det(w,x)^2}{|x+sw|^3}\ge 0.$$
In particular, the function $\phi:s\longmapsto Q_2Q_1(s)$ is strictly convex as for every $s$ sufficiently small 
$$\phi''(s) = \frac{d^2}{ds^2}|Q_1+s\Vec{u}-Q_2| = \frac{\det(\overrightarrow{Q_2Q_1},\Vec{u})^2}{|\overrightarrow{Q_2Q_1}+s\Vec{u}|^3}>0,$$
where the strict positivity follows from the fact that $\Vec{u}$ and $Q_1-Q_2$ are not collinear. 

This allows to conclude that $\Psi$ is strictly convex in the neighborhood of $0$. Therefore, there exists a direction ($s>0$ or $s<0$) for which 
$$P(K_s)>P(K)\ \ \text{with}\ \ |K_s|=|K|,$$
for sufficiently small $s$. Therefore, the quadrilateral $K$ cannot be a solution of \eqref{problem}, which means that a solution sould be a triangle.

\vspace{2mm}

\textbf{Step 4: The optimal triangle $K$ shares at least one vertex with $T$}\vspace{2mm}

Assume by contradiction that an optimal  triangle $K$ does not share any vertex with $T$, and denote by $Q_1$, $Q_2$ and $Q_3$ the vertices of the triangle $K$. All of these vertices are located in the interior of sides of $T$. We then have two cases: 
\begin{itemize}
    \item either each side of $T$ contains exactly one vertex of $K$,
    \item  or there exists a side that contains exactly two vertices.
\end{itemize}
 The latter case can be easily excluded by considering a subset $\Tilde{K}\subset T$ which is the image of $K$ by a (small) rotation. At least one of the vertices of $\Tilde{K}$ will be located inside $T$, contradicting the optimality of $\Tilde{K}$ and therefore the one of $K$. In what follows, we consider the first case.     

We denote by $A$, $B$ and $C$ the vertices of $T$ and assume first that there exists $\alpha,\beta,\gamma\in (0,1)$ such that 
\[
Q_1=(1-\alpha)B+\alpha C,\qquad
Q_2=(1-\beta)C+\beta A,\qquad
Q_3=(1-\gamma)A+\gamma B.
\]
Set
$\delta:=\beta+\gamma-1$ and $
R:=\sqrt{(1-\delta)^2+4\alpha\delta}$ and consider  
\[
\lambda:=\frac{2\delta}{1-\delta+R}
\qquad \text{and}\qquad
\mu:=-\frac{2\delta}{1+\delta+R}.
\]
These coefficients are well-defined since $\delta\in(-1,1)$ and $R\ge 0$. 

Consider the perturbed vertices
\[
\begin{aligned}
Q_1(s)&:=Q_1+s\overrightarrow{BC},\\
Q_2(s)&:=Q_2+s\lambda\overrightarrow{AB},\\
Q_3(s)&:=Q_3+s\mu\overrightarrow{CA},
\end{aligned}
\]
and let $K_s$ the triangle of vertices $Q_1(s)$, $Q_2(s)$ and $Q_3(s)$.
For sufficiently small $|s|$, a direct determinant computation gives
\[
|K_s|-|K|
=|T|\bigl(\delta-\alpha\lambda+(1-\alpha)\mu\bigr)s-|T|\bigl(\lambda+\mu+\lambda\mu\bigr)s^2.
\]

The above choice of $\lambda$ and $\mu$ satisfies
\begin{equation*}
   \delta-\alpha\lambda+(1-\alpha)\mu=0
\qquad \text{and}\qquad
\lambda+\mu+\lambda\mu=0. 
\end{equation*}

Consequently, $|K_s|=|K|$ and $K_s\subset T$, for all sufficiently small $|s|$.

Let us now show that the function $s\mapsto P(K_s)$ is strictly convex in a neighborhood of $0$.

Set
\[
x:=\overrightarrow{Q_1Q_2},\qquad
y:=\overrightarrow{Q_1Q_3},\qquad
W:=(1+\lambda)\overrightarrow{AB}-\overrightarrow{AC}.
\]
Since $T$ is nondegenerate, the vectors $B-A$ and $C-A$ are
linearly independent. Consequently,
$W\neq0$. Moreover, $Q_1,Q_2,Q_3$ are noncollinear, so $\det(x,y)\neq0$.

Using $\lambda+\mu+\lambda\mu=0$, we can write
\[
\begin{aligned}
\overrightarrow{Q_1(s)Q_2(s)}&=x+sW,\\
\overrightarrow{Q_2(s)Q_3(s)}&=y-x+\mu sW,\\
\overrightarrow{Q_1(s)Q_3(s)}&=y+(1+\mu)sW.
\end{aligned}
\]
Therefore, for all sufficiently small $|s|$, we have 
\[
P(K_s)
=
| x+sW |
+| y-x+\mu sW |
+ | y+(1+\mu)sW |.
\]

The second order derivative of $s\longmapsto P(K_s)$ in the neighborhood of $0$ is given by 
\[
\frac{d^2}{ds^2}P(K_s)
=
\frac{\det(x,W)^2}{| x+sW|^3}+\mu^2
\frac{\det(y-x,W)^2}{| y-x+\mu sW|^3}+(1+\mu)^2
\frac{\det(y,W)^2}{| y+(1+\mu)sW|^3}.
\]

Since $x$ and $y$ are not collinear (as the triangle $K$ is supposed to be nondegenerate), we deduce that at least one of them is not collinear to $W$, which means that at least one of the terms $\det(x,W)^2$ and $\det(y,W)^2$ is strictly positive. Therefore, the function $s\longmapsto P(K_s)$ is strictly convex in the neighborhood of $0$ and therefore there exists a direction ($s>0$ or $s<0$) such that for small $|s|$, we have 
$$P(K_s)>P(K)\qquad \text{with}\qquad |K_s|=|K|.$$
This contradicts the maximality of $K$ as a solution of \eqref{problem}.  Therefore, an optimal triangle shares at least one vertex with $T$. \newpage

\textbf{Step 5: The optimal triangle shares two vertices with $T$  }\vspace{2mm}

Without loss of generality, we can assume that the triangle $K=Q_1Q_2Q_3$ shares exactly one vertex with $T$ and that vertex is $A=Q_1$. We then have to exclude the following three cases: 
\begin{itemize}
    \item $Q_2$ is located in the interior the side $[AB]$ and $Q_3$ is in the interior of the side $[AC]$. 
    \item $Q_2$ is located in the interior the side $[AB]$ and $Q_3$ is in the interior of the side $[BC]$.
    \item $Q_2$ and $Q_3$ are located in the interior of the side $[BC]$.
\end{itemize}

The first case can be easily excluded, by considering $K_s$ a translation of $K$ by the vector $s(B-A)$, with sufficiently small $s>0$. One of the vertices of $K_s$ is included in the interior of $T$, which excludes the optimality of $K_s$ and therefore also the one of $K$ by following \textbf{Step 1}. 

The second case can also be excluded by considering $K_\theta$ the image of $K$ with a rotation of center $A=Q_1$ and angle $\theta$. For sufficiently small $|\theta|$ and the adequate sign of $\theta$, the triangle $K_\theta$ is included in the container $T$ and has two vertices located in the interior of $T$, which excludes its optimality and therefore also the one of $K$.   

\begin{figure}[h!]
    \centering
    \includegraphics[width=0.5\linewidth]{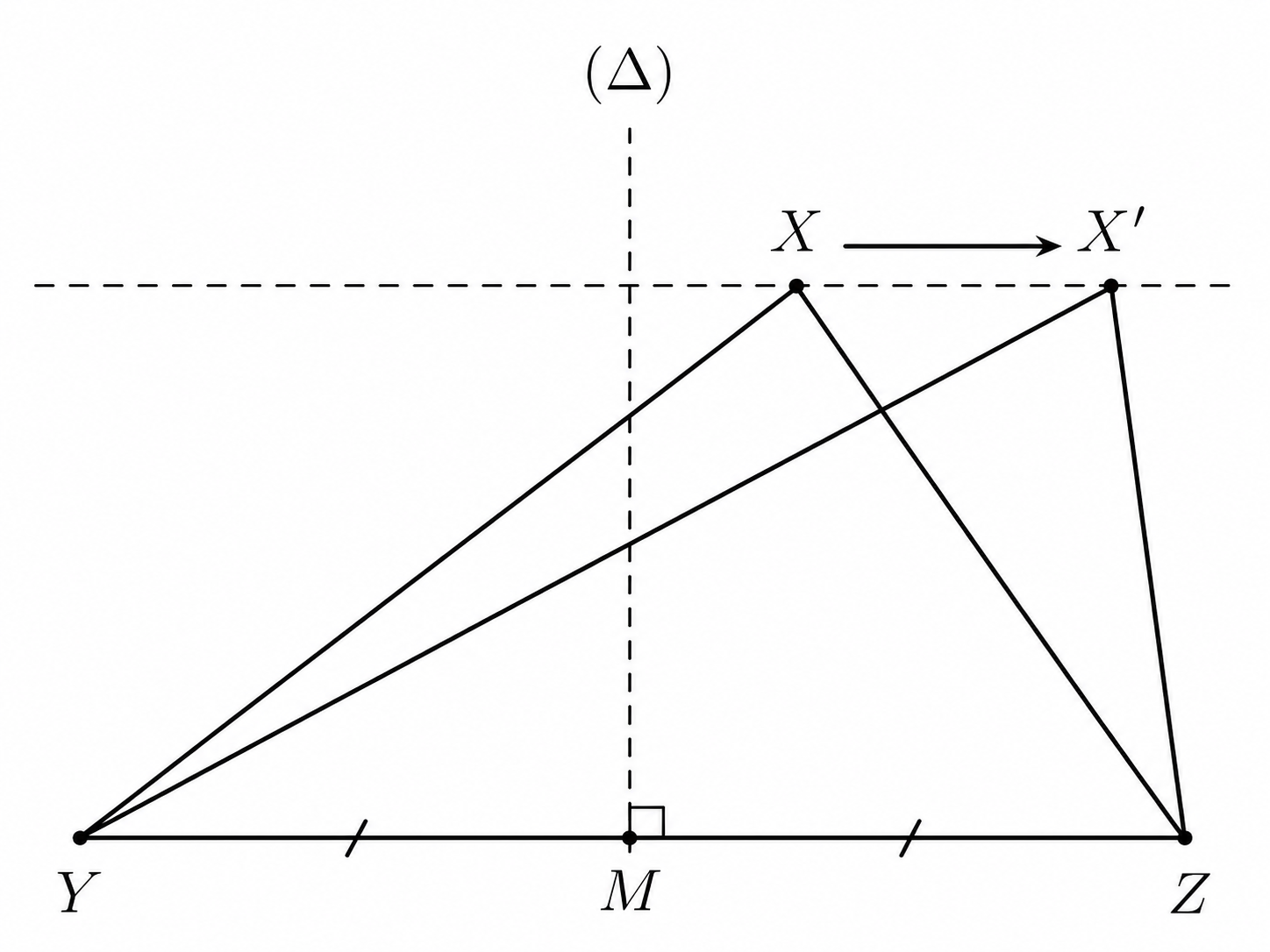}
    \caption{Parallel chord movements increase (or decrease) the perimeter while keeping the area fixed.}
    \label{fig:parallel}
\end{figure}

To exclude the third case, we will make use of a special kind of transformation called \textit{parallel chord movement}. The idea of such transformation is the following: consider $XYZ$ a triangle of vertices $X$, $Y$ and $Z$ and $(\Delta)$ the perpendicular bisector of the segment $[YZ]$. Assume that $X'$ is obtained by translating $X$ along the line containing $X$ and parallel to the line $(YZ)$, and consider $X'YZ$ the triangle of vertices $X'$, $Y$ and $Z$, see Figure \ref{fig:parallel}. It is clear that $|XYZ|=|X'YZ|$ (because they share the same altitude with respect to the common side $[YZ]$) and that if $d(X,(\Delta))< d(X',(\Delta))$, then, $P(X'YZ)>P(XYZ)$. 

Therefore, considering $Q_2(s):= Q_2 + s \overrightarrow{BC}$ and $Q_3(s):= Q_3 + s \overrightarrow{BC}$, for small $|s|$, we obtain a perturbation $K_s:=Q_1Q_2(s)Q_3(s)$ of the triangle $K$ that is included in the contained $T$, has the same area as $K$, and is the image of $K$ with a translation and a parallel chord movement. Then one just has to choose the suitable sign of $s$ that gives $P(K_s)>P(K)$, which contradicts the optimality of $K$. 

We then have showed that the optimal triangle shares exactly two vertices with the container $T$.\vspace{2mm}

\textbf{Step 6: The best triangle }\vspace{2mm}

Let us denote $a:= BC$, $b:=AC$, $c:=AB$, $\alpha := \widehat{BAC}$, $\beta := \widehat{ABC}$ and $\gamma := \widehat{ACB}$ and assume that $a\ge b\ge c$, which yields $\alpha\ge \beta\ge \gamma$. We have brought the problem to the comparison of the perimeters of 6 triangles of the same area $\eta |T|$ each of them sharing exactly two vertices with the triangle $T$ as shown in Figure \ref{fig:six_cases}. For every $(X, \theta)\in \{A,B,C\}\times \{\alpha,\beta,\gamma\}$, we assign a triangle $K_{X,\theta}$, where $X$ stands for the vertex that is not common between $K_{X,\theta}$ and $T$ and $\theta$ to the common angle between $K$ and $T$.  

\begin{table}[h!]
    \centering
    \begin{tabular}{|c|c|}
        \hline
        \includegraphics[width=0.35\linewidth]{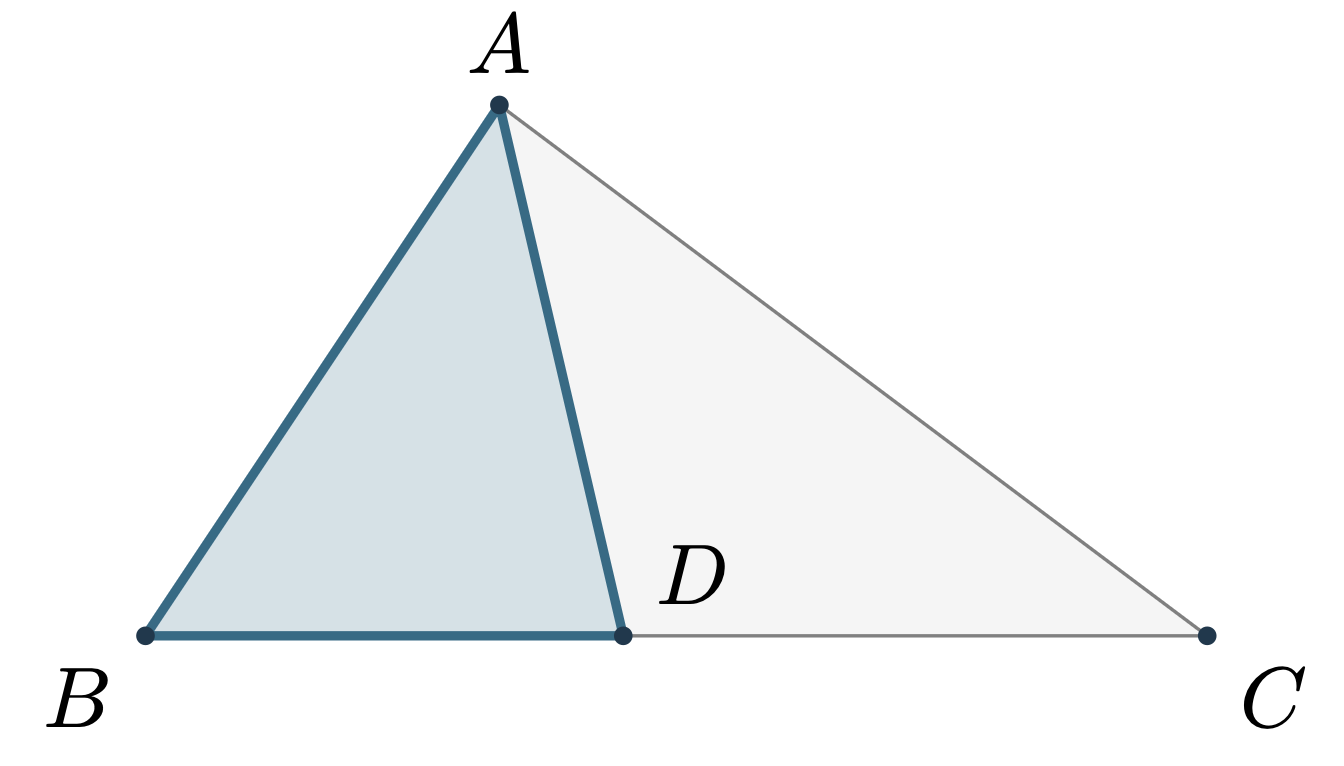} &
        \includegraphics[width=0.35\linewidth]{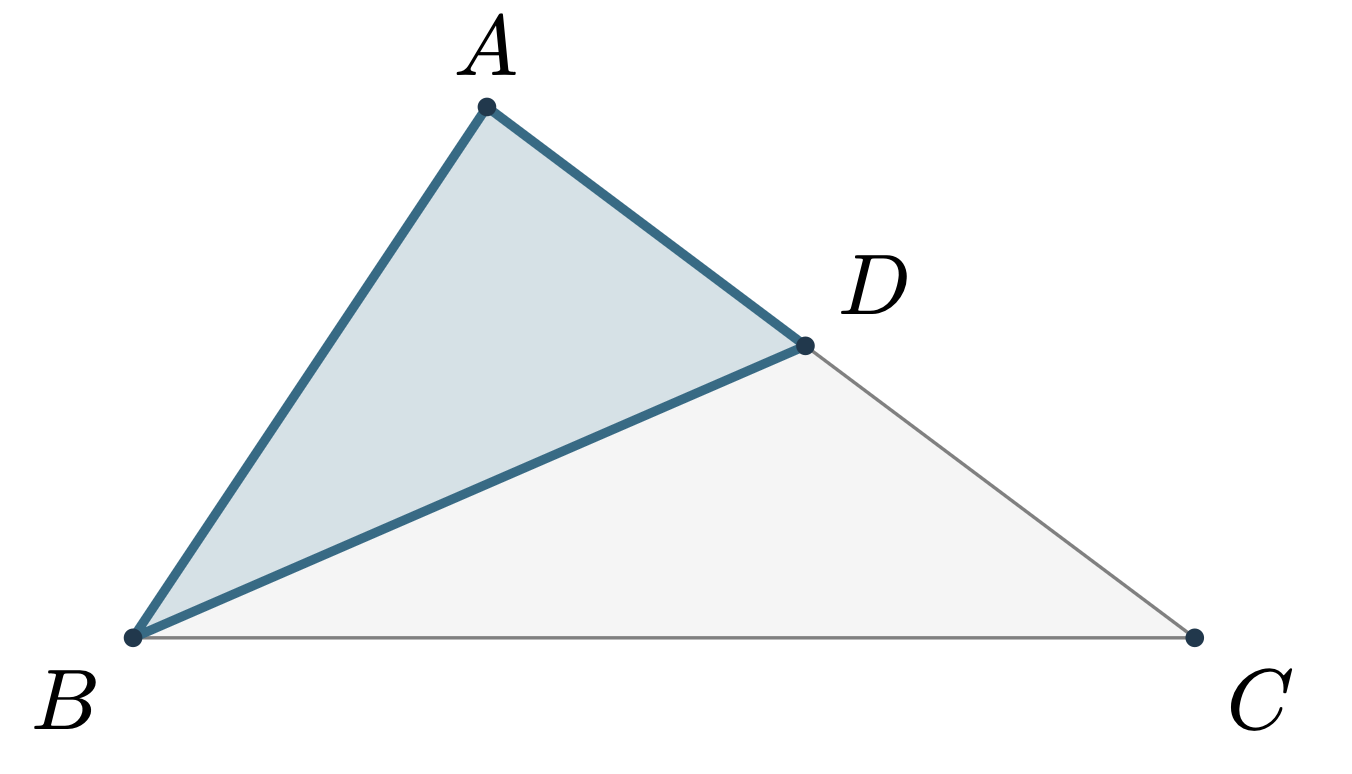} \\
        \hline
        Triangle $K_{C,\beta}$ & Triangle $K_{C,\alpha}$ \\
        \hline
        \includegraphics[width=0.35\linewidth]{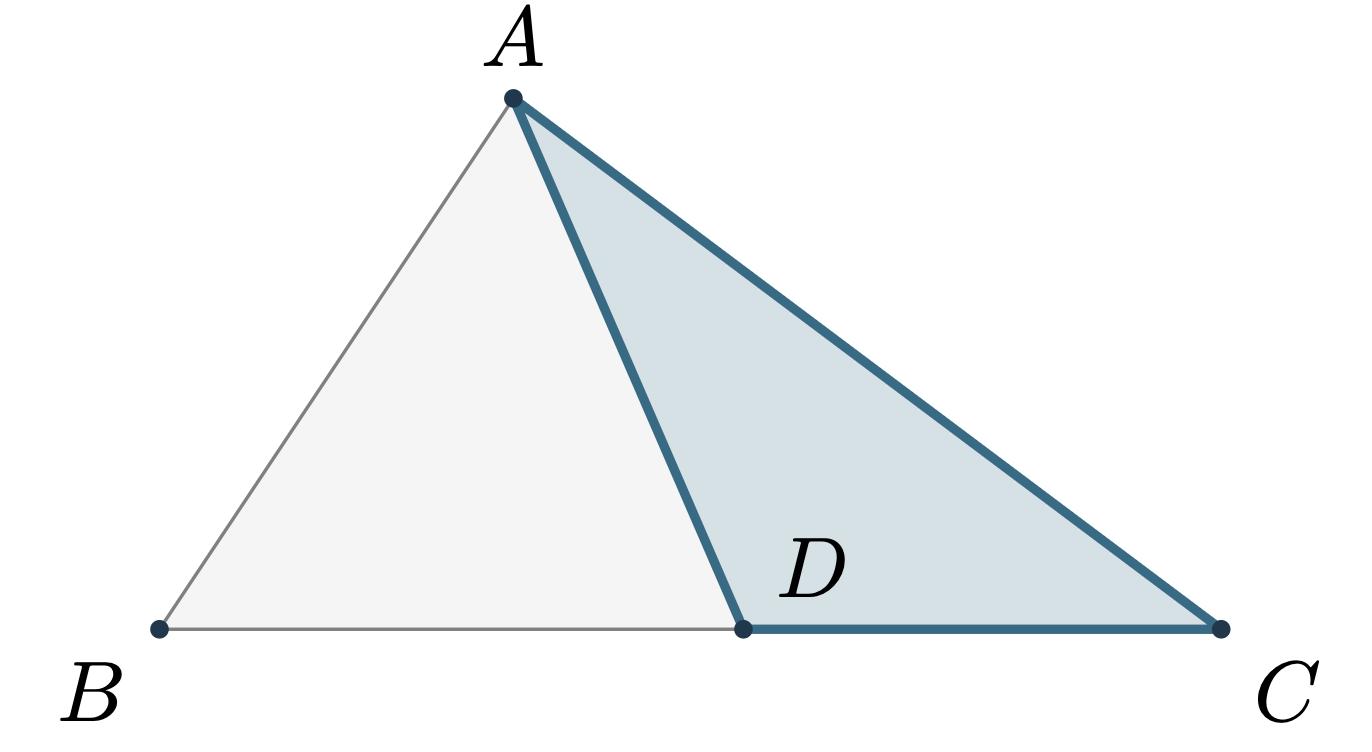} &
        \includegraphics[width=0.35\linewidth]{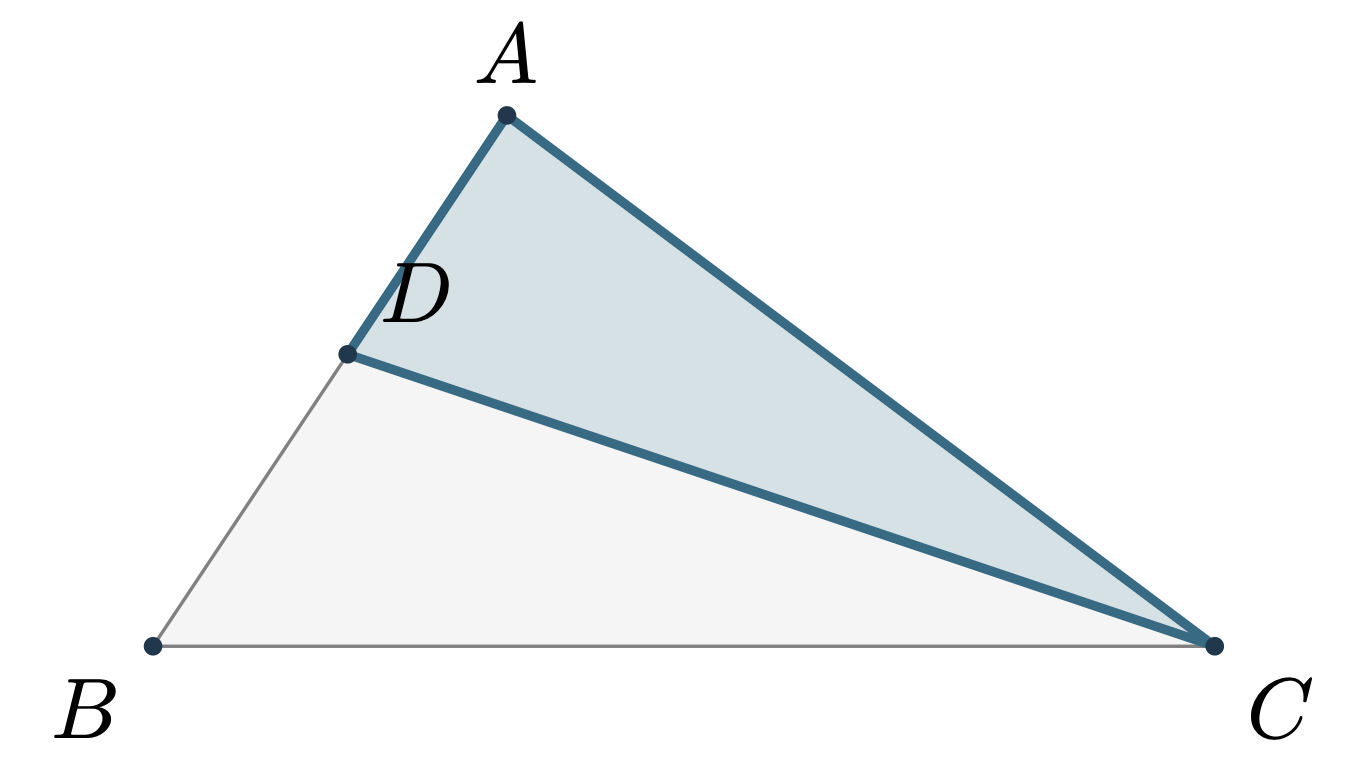} \\
        \hline
        Triangle $K_{B,\gamma}$ & Triangle $K_{B,\alpha}$ \\
        \hline
        \includegraphics[width=0.35\linewidth]{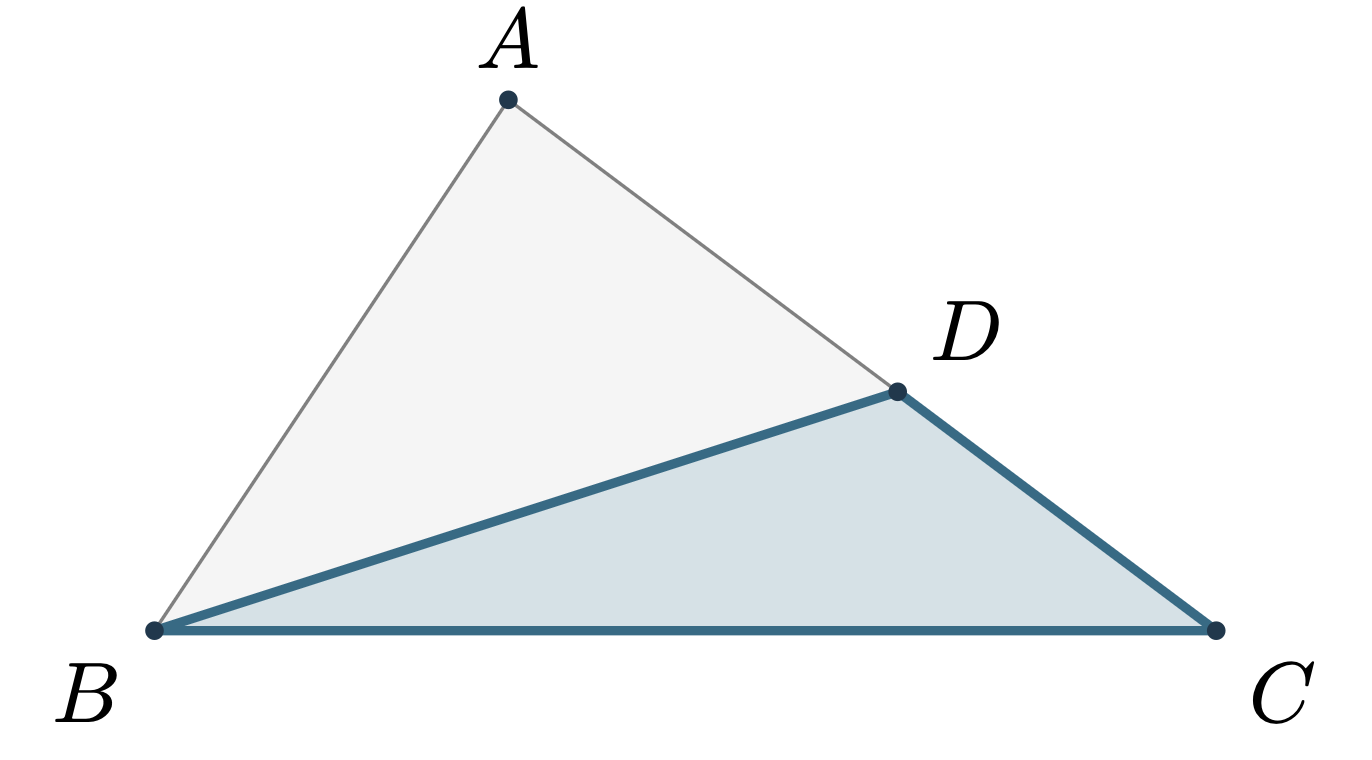} &
        \includegraphics[width=0.35\linewidth]{figures/6.png} \\
        \hline
         Triangle $K_{A,\gamma}$ & Triangle $K_{A,\beta}$  \\
        \hline
    \end{tabular}
    \caption{The remaining 6 triangles all of the same area $\eta |T|$.}
    \label{fig:six_cases}
\end{table}

Let us now show that the best configuration, i.e., with the maximal perimeter, is the configuration $(A,\beta)$. To do so, we will make use of \textit{parallel chord movements} introduced in \textbf{Step 5}.

By using simple reflection and parallel chord movement arguments, we can show that $P(K_{A,\beta})\ge P(K_{A,\gamma})$, with equality if and only if $\beta=\gamma$. Indeed, consider $D_\beta$ and $D_\gamma$ the vertices of $K_{A,\beta}$ and $K_{A,\gamma}$ respectively, $(\Delta)$ the perpendicular bisector of the segment $[BC]$ and $D_\gamma'$ the symmetric of $D_\gamma$ with respect to the line $(\Delta)$ and $K'_{A,\gamma}$ the triangle of vertices $D_\gamma'$, $B$ and $C$. We have $|K_{A,\gamma}|=|K_{A,\beta}|=|K'_{A,\gamma}|$ and since $d(D_\beta,(\Delta))\ge d(D'_\gamma,(\Delta))$ (with equality if and only if $\beta=\delta$), we deduce that $$P(K_{A,\gamma})=P(K'_{A,\gamma})\leq P(K_{A,\beta}),$$ with equality if and only if $\beta=\delta$. 

Using similar arguments, we prove that
\begin{itemize}
    \item $P(K_{B,\alpha})\ge P(K_{B,\gamma})$, with equality if and only if $\alpha = \gamma$.
    \item $P(K_{C,\alpha})\ge P(K_{C,\beta})$, with equality if and only if $\alpha = \beta$.
\end{itemize}

It remains then to compare $P(K_{A,\beta})$, $P(K_{C,\alpha})$ and $P(K_{B,\alpha})$. We denote by $D_\beta$ the third vertex of the triangle $K_{A,\beta}$ (along with the vertices $B$ and $C$) and $D_\alpha$ the third vertex of the triangle $K_{B,\alpha}$ (along with the vertices $A$ and $B$). We remark that the triangles $K_{A,\beta}$ and $K_{B,\alpha}$ share the same altitude (with respect to the common line $(AB)$), therefore, since they have the same area, we deduce that $AD_\alpha = BD_\beta$. We then define $\widetilde{K}$ as the triangle obtained by translating $K_{A,\beta}$ by the vector $\overrightarrow{D_\beta A}$ and $\Tilde{C}:= C+ \overrightarrow{D_\beta A}$. By noting that $\widetilde{K}$ is the image of $K_{B,\alpha}$ by a parallel chord movement such that $d(\Tilde{C},(\Delta))\ge d(C,(\Delta))$, where $(\Delta)$ is the orthogonal bisector of the segment $[AD_\alpha]$, and with equality if and only if $\beta=\alpha$, we deduce that 
$$P(K_{A,\beta})=P(\widetilde{K})\ge P(K_{B,\alpha}),$$
with equality if and only if $\alpha=\beta$ (in which case $\widetilde{K}$ and $K_{B,\alpha}$ will be either identical or symmetric with respect to $(\Delta)$). 

With the exact similar arguments, we can show that 
$$P(K_{A,\gamma})\ge P(K_{C,\alpha}),$$
with equality if and only if $\alpha=\gamma$. 

We finally conclude that the triangle $K_{A,\beta}$ is an optimal shape for \eqref{problem} and that it is unique if and only if $\alpha>\max(\beta, \gamma)$, which concludes the proof.

\medskip\noindent



\bibliographystyle{abbrv}
\bibliography{biblio.bib}

\end{document}